\documentclass[12pt]{article}   

\usepackage{amssymb}     
\usepackage{latexsym}       
\newcommand{\bq}{\begin{quote}}
\newcommand{\eq}{\end{quote}}
 
\newcommand{\lood}{^\perp}
\newcommand{\acht}{\wp}

\newcommand{\pijl}{\longrightarrow}

\newcommand{\M}{\mid}
\newtheorem{Th}{Theorem}  
\newtheorem{ax}{Axiom}  
\newtheorem{lm}{Lemma} 
\newtheorem{df}{Definition}    
\newtheorem{pr}{Proposition} 
\newtheorem{cl}{Corollary}  
\newtheorem{re}{Remark}    
\newtheorem{as}{Assumption}  
\newtheorem{wg}{Wild Guess}
\newtheorem{ex}{Example}
\newcommand{\bth}{\begin{Th}\hspace{-5pt}{\bf .} \ } 
\newcommand{\Eth}{\end{Th}}
\newcommand{\bax}{\begin{ax}\hspace{-5pt}{\bf .} \ } 
\newcommand{\eax}{\end{ax}}
\newcommand{\blm}{\begin{lm}\hspace{-5pt}{\bf .} \ }
\newcommand{\elm}{\end{lm}}
\newcommand{\bdf}{\begin{df}\hspace{-5pt}{\bf .} \ }   
\newcommand{\edf}{\end{df}} 
\newcommand{\bpr}{\begin{pr}\hspace{-5pt}{\bf .} \ } 
\newcommand{\epr}{\end{pr}}
\newcommand{\bcl}{\begin{cl}\hspace{-5pt}{\bf .} \ } 
\newcommand{\ecl}{\end{cl}}
\newcommand{\bre}{\begin{re}\hspace{-5pt}{\bf .} \ }
\newcommand{\ere}{\end{re}}
\newcommand{\bas}{\begin{as}\hspace{-5pt}{\bf .} \ }
\newcommand{\eas}{\end{as}}
\newcommand{\bwg}{\begin{wg}\hspace{-5pt}{\bf .} \ }
\newcommand{\ewg}{\end{wg}}
\newcommand{\bex}{\begin{ex}\hspace{-5pt}{\bf .} \ }  
\newcommand{\eex}{\end{ex}}

\newcommand{\bit}{\begin{itemize}}
\newcommand{\eit}{\end{itemize}\par\noindent}
\newcommand{\ben}{\begin{enumerate}}
\newcommand{\een}{\end{enumerate}\par\noindent}
\newcommand{\beq}{\begin{equation}}
\newcommand{\eeq}{\end{equation}\par\noindent}
\newcommand{\beqa}{\begin{eqnarray*}}
\newcommand{\eeqa}{\end{eqnarray*}\par\noindent}
\newcommand{\beqn}{\begin{eqnarray}}  
\newcommand{\eeqn}{\end{eqnarray}\par\noindent}

\ifx\optionkeymacros\undefined\else\endinput\fi
\catcode`\ç=\active\defç{\c c}        
\catcode`\¾=\active\def¾{\leq}        
\catcode`\"=\active\def"{\geq}        
\catcode`\³=\active\def³{``}          
\catcode`\Œ=\active\defŒ{`}           
\catcode`\¹=\active\def¹{'}           
\catcode`\Ç=\active\defÇ{\c C}        
\catcode`\'=\active\def'{{\OE}}       
\catcode`\°=\active\def°{\accent'27}  
\catcode`\²=\active\def²{''}          
\catcode`\ä=\active\defä{\"a}        
\catcode`\ë=\active\defë{\"e}        
\catcode`\ï=\active\defï{\"{\i}}     
\catcode`\ö=\active\defö{\"o}        
\catcode`\ü=\active\defü{\"u}        
\catcode`\ÿ=\active\defÿ{\"y}        
\catcode`\Ä=\active\defÄ{\"A}        
\catcode`\Ö=\active\defÖ{\"O}        
\catcode`\Ü=\active\defÜ{\"U}        
\catcode`\á=\active\defá{\'a}        
\catcode`\é=\active\defé{\'e}        
\catcode`\í=\active\defí{\'{\i}}     
\catcode`\ó=\active\defó{\'o}        
\catcode`\ú=\active\defú{\'u}        
\catcode`\É=\active\defÉ{\'E}        
\catcode`\à=\active\defà{\`a}        
\catcode`\è=\active\defè{\`e}        
\catcode`\ì=\active\defì{\`{\i}}     
\catcode`\ò=\active\defò{\`o}        
\catcode`\ù=\active\defù{\`u}        
\catcode`\À=\active\defÀ{\`A}        
\catcode`\ã=\active\defã{\~a}        
\catcode`\ñ=\active\defñ{\~n}        
\catcode`\õ=\active\defõ{\~o}        
\catcode`\Ã=\active\defÃ{\~A}        
\catcode`\Ñ=\active\defÑ{\~N}        
\catcode`\Õ=\active\defÕ{\~O}        
\catcode`\â=\active\defâ{\^a}        
\catcode`\ê=\active\defê{\^e}        
\catcode`\î=\active\defî{\^{\i}}     
\catcode`\ô=\active\defô{\^o}        
\catcode`\û=\active\defû{\^u}        
\catcode`\æ=\active\defæ{{\ae}}      
\let\optionkeymacros\null

\begin{document}  
\noindent\centerline{\LARGE{\bf Logic of Dynamics $\&$ Dynamics of Logic;}} 
 
\smallskip\noindent\centerline{\LARGE{\bf Some Paradigm Examples}}

\bigskip\par\noindent
\medskip\par\noindent
\centerline{\normalsize{Bob Coecke}}
\par\medskip\noindent
\centerline{\footnotesize{University of Oxford, Oxford University Computing Laboratory\,,}}\vspace{-1mm}
\par\noindent 
\centerline{\footnotesize{Wolfson Building, Parks Road, Oxford, OX1 3QD, UK\,;}}\vspace{-1mm}
\par\noindent 
\centerline{\footnotesize{e-mail: bocoecke@vub.ac.be \& coecke@comlab.ox.ac.uk\,.}}\vspace{-1mm} 
 
\par\bigskip\noindent  
\centerline{{David J. Moore}}
\par\medskip\noindent
\centerline{\footnotesize{University of Canterbury, Department of Physics and   
Astronomy,}}\vspace{-1mm}
\par\noindent 
\centerline{\footnotesize{Private Bag 4800, Christchurch\,;
New Zealand.}}\vspace{-1mm}
\par\noindent 
\centerline{\footnotesize{e-mail: closcat@hotmail.com\,.}}
  
\par\bigskip\noindent 
\centerline{{Sonja Smets}}
\par\medskip\noindent
\centerline{\footnotesize{Free University of
Brussels (VUB), Department of Philosophy,}}\vspace{-1mm}
\par\noindent 
\centerline{\footnotesize{Pleinlaan 2, B-1050
Brussels, Belgium\,;}}\vspace{-1mm} 
\par\noindent 
\centerline{\footnotesize{e-mail: sonsmets@vub.ac.be\,.}} 
\par\bigskip\noindent
%

\def\retro{\hbox{$\small\, \circ\hspace{-1mm} -\hspace{-1mm}\ $}}  

\par\bigskip\noindent
\begin{abstract} 
\par\noindent
The development of ``operational quantum logic" points out
that classical boolean structures are too rigid to describe
the actual and potential properties of quantum systems.
Operational quantum logic bears upon basic axioms which are
motivated by empirical facts and as such supports the
dynamic shift from classical to non-classical logic
resulting into a dynamics of logic.

\smallskip
On the other hand, an intuitionistic perspective on operational
quantum logic, guides us in the direction of incorporating dynamics logically by reconsidering the
primitive
propositions required to describe the behavior of a quantum system, in
particular in view of the emergent disjunctivity due to the non-determinism of
quantum measurements.

\smallskip
A further elaboration on ``intuitionistic quantum logic" emerges into a ``dynamic
operational quantum logic", which allows us to express dynamic
reasoning in the sense that we can capture how actual
properties propagate, including their temporal causal
structure.  It is in this sense that passing from static
operational quantum logic to dynamic operational quantum
logic results in a true logic of dynamics that provides a
unified logical description of systems which evolve or which
are submitted to measurements.  This setting reveals that even
static operational quantum logic bears a hidden dynamic ingredient in
terms
of what is called ``the orthomodularity" of the lattice-structure.

\smallskip
Focusing on the
quantale semantics for dynamic operational quantum logic, we
delineate some points of difference with the existing
quantale semantics for (non)-commutative linear logic.
Linear logic is here to be conceived of as a resource-sensitive
logic capable of dealing with actions or in other words, it
is a logic of dynamics.

\smallskip
We take this opportunity to
dedicate this paper to Constantin Piron at the occasion of
his retirement.
\end{abstract} 

\bigskip\noindent
{\bf 1. INTRODUCTION}
     
\medskip\noindent
As a starting point for our discussion on the dynamics of logic we quote G. Birkhoff and J. von Neumann,
confronting the then ongoing tendencies towards intuitionistic logic with their observation of the
``logical'' structure encoded in the lattice of closed subspaces of a Hilbert space, the
``semantics'' of quantum theory (Birkhoff and von Neumann 1936)\,:
\begin{quote}
``The models for propositional calculi [of physically significant
statements in quantum
mechanics] are also interesting from the standpoint of pure logic.
Their nature is determined by quasi-physical and technical reasoning,
different from the
introspective and philosophical considerations which have to guide
logicians hitherto [\,...\,]
whereas logicians  have usually assumed that [the
orthocomplementation] properties
L71-L73 of negation were the ones least able to withstand a critical
analysis, the study of
mechanics points to the distributive identities L6 as the weakest
link in the algebra
of logic.'' (p.839)
\end{quote}   
They point at a fundamental 
difference between Heyting algebras (the semantics of
intuitionistic propositional logic) and orthomodular lattices (the
``usual'' semantics of
quantum logic) when  viewed as generalizations of Boolean algebras (the semantics of 
classical propositional logic). A new intuitionistic perspective on operational quantum logic (see below) provides a way
of blending these seemingly contradicting directions in which logic propagated during the previous century (Coecke
2002)\,. In this paper we focus on two new logical structures, namely (intuitionistic) linear logic (Girard 1987, Abrusci
1990) which emerged from the traditional branch of logic, and dynamic operational quantum logic (Coecke (nd), Coecke
and Smets 2001), emerging from an elaboration on the above mentioned blend. We also briefly consider the ``general
dynamic logic'' proposed in van
Benthem (1994)\,. We do mention epistemic action logic (e.g., Baltag 1999) and computation and information flow (e.g.,
Abramsky 1993, Milner 1999) as other examples of dynamic aspects in logic, which we unfortunately will not be able to
consider in this paper. 
We also won't discuss the ``geometry of interaction'' paradigm which provides a (promising) different perspective on
linear logic (Abramsky and Jagadeesan 1994), but is still in full development. Concretely we start in section 2 with
an outline of static operational quantum logic. In section 3 we survey dynamic operational quantum logic,
and demonstrate that the emerging structures do not fit in the logical system proposed in van
Benthem (1994)\,. In section 4 we compare dynamic operational quantum logic with linear logic while we focus
on formal and methodological differences.  

\bigskip\noindent 
{\bf 2. STATIC ``OLD-STYLE'' OPERATIONAL QUANTUM LOGIC}
     
\medskip\noindent
The Geneva School approach to the logical foundations of physics originated with the work in
Jauch and Piron (1963), Piron (1964), Jauch (1968), Jauch and Piron (1969), Piron (1976) and Aerts
(1981) as an incarnation of the part of the research domain ``foundations of physics'' that is
nowadays called ``Operational quantum logic" (OQL)$^1$ --- see Coecke et al (2000) for a recent overview of its
general aspects. With respect to the intuitionistic perspective and the dynamic aspects which we put forward
below, we will further refer in this paper to OQL as ``sOQL", emphasizing its static nature.
More concretely, sOQL as a theory aims to characterize physical systems, ranging from classical to
quantum, by means of their actual and potential properties, in particular by taking an
ontological rather than an empirical perspective, but, still providing a truly operational
alternative to the standard approaches on the logical status of quantum theory.$^2$ Since Moore (1999),
Coecke et al ((nd)a,b) and Smets (2001) provide recent and detailed discussions of sOQL, we
intend in this paper to give only a
concise overview, focusing on its basic concepts and underlying epistemology.  
We are aware of the
strong conceptual restrictions imposed by the rigid foundation of sOQL,
necessary in order to obtain a framework with fully
well-defined primitive notions. Clearly, in view of still existing
conceptual incompatibilities at the foundational level between quantum
theory
and relativistic space-time and field theoretic considerations,
the development of an essentially ``towards dynamics directed"-formalism for quantum
logicality should range beyond the rigid concepts of sOQL.  Therefore, we conceive of the notions
inherited
from sOQL as a stepping-stone for further development rather than something
necessarily ``to carry all the way''. This is the reason why we lowercase `o' in our notation DoQL
referring to ``dynamic operational quantum logic" and IoQL, when referring to  ``intuitionistic operational
quantum logic". In particular the sOQL assumptions (see below) of ``a particular physical system which is
considered as distinct from its surroundings", ``the specification of states as
a (pre-defined) set" and ``the a priori specification of a particular physical
system itself" need too be reconsidered when crossing the edges of sOQL  (which itself was designed to clarify the
structural description of quantum systems and justify an ontological perspective for non-relativistic
quantum theory).
 
\smallskip
First we want to clarify that the operationalism which forms the core of sOQL points to a pragmatic attitude
and not to any specific doctrine one can encounter in last century's philosophy of science (Coecke et al
(nd)a).  Linked to the fact that we defend the position of critical scientific realism in relation of sOQL,
operationality points to the underlying assumption that with every property of a physical system we can
associate experimental procedures that can be performed on the system,  each such experimental procedure
including the specification of a well-defined positive result for which certainty is exactly guaranteed by
the actual existence of the corresponding property.  In particular, while on the epistemological level our
knowledge of what exists is based on what we could measure or observe, on the ontological level 
``physical" properties have an extension in reality and are not reducible to sets of procedures.$^3$
Focusing on the stance of critical scientific realism, we first adopt an ontological realistic position and
a correspondence theory of truth.$^4$  Though, contrary to naive realists we do adopt the thesis of
fallibilism by which truth in relation to scientific theories has to be pursued but can only be approached. 
As such we agree with Niiniluoto (1999) that scientific progress can be characterized in terms of increasing
truthlikeness.  Furthermore, we believe that reality can indeed be captured in conceptual frameworks, though
contrary to Niiniluoto we like to link this position to a ``weak" form of conceptual idealism. 
Concentrating for instance on Rescher's conceptual idealism as presented in Rescher (1973, 1987
\S 11) and revised in Rescher (1995 \S 8), this position maintains that a description of physical reality
involves reference to mental operations, it doesn't deny ontological realism and also adopts the thesis of
fallibilism. This position is opposed to an ontological idealism in which the mind produces the ``real"
objects. As one might expect, we want to stress that we are not inclined to adopt Rescher's strict Kantian
distinction between reality out there and reality as we perceive it, though we are attracted by the idea of
capturing  ontological reality in mind-correlative conceptual frameworks.  Once we succeed in giving such a
description of reality, it approaches according to us the ontological world close enough to omit a Kantian
distinction between the realms of noumena and phenomena.  Why it is of importance for us to reconcile
critical realism with a weak form of conceptual idealism becomes clear when we focus on
sOQL.  Firstly, we cannot escape the fact that in our theory we focus on parts of reality considered as well-defined
and distinct which we can then characterize as
physical systems.  Secondly we have to identify the properties of those systems which is a mind-involving
activity.   To be more explicit, a physicist can believe that a system ontologically has
specific actual and potential properties, but to give the right characterization of the
physical system he has to consider the definite experimental projects which can test those
properties so that, depending on the results he would obtain, he can be reinforced in his
beliefs or has to revise them.$^5$ How sOQL formally is built up, using the notions of actual properties,
potential properties and definite experimental projects, will be clarified below.  To
recapitulate we finish this paragraph by stressing that we focus in our scientific activity
on parts of the external world, reality is mind-independent, though the process of its
description ``involves" some mind-dependent characterizations.

\smallskip
Let us introduce the primitive concrete notions on which sOQL relies, explicitly following
Coecke et al ((nd)a):
\bit
\item
We take a particular physical system to be a part of the ostensively external world which is considered as distinct from its
surroundings --- see Moore (1999) for a discussion on this matter;  
\item
A singular realization of the given particular physical system is a conceivable manner of being of that system within a
circumscribed experiental context;
\item
States ${\cal E} \in \Sigma$ of a given particular physical system are construed as abstract names encoding its possible
singular realizations;
\item
A definite experimental project $\alpha \in {\cal Q}$ on the given particular physical system is a real experimental
procedure which may be effectuated on that system where we have defined in advance what would be the positive response
should it be performed.
\item
Properties $a \in {\cal L}$ of a given particular physical system are construed as 
{\it candidate}  elements of
reality  corresponding to the definite experimental projects defined for
that system.  We as such obtain a mapping of definite experimental projects ${\cal Q}$ on
properties
${\cal L}$\,.  
\eit
The notion of an element of reality was first introduced in Einstein et
al (1935) as follows: 
\begin{quote}
``If, without in any way disturbing a system, we can predict with certainty [...] the value of a
physical quantity, then there exists an element of physical reality corresponding to this physical quantity"
(p.777). 
\end{quote}
While this formulation is explicitly the starting point for Piron's early work, we insist that 
Piron's operational concept of element of reality is both more precise and allows theoretical deduction, being based
on an empirically accessible notion of counterfactual performance rather than a metaphorical notion of
non-perturbation.  The basic ingredient that we inherit
from this setting is the agreement that actual properties exist.  Before we give a more rigid characterization of how
an actual property is conceived within sOQL, we want to remark that Piron adopted the Aristotelian
concepts of actuality and potentiality and placed them in a new framework (see e.g. Piron (1983)). For an
analysis of how these notions, used within sOQL, are still related to the old Aristotelian ones, we refer to
Smets (2001).   Let us just briefly mention here that with regard to actuality, an actual property is within
sOQL conceived as an attribute which ``exists"; it is some realization in reality or in other words;
an element of reality.  A potential property on the other hand, does not exist in the same way as an actual
one, it is conceived merely as a capability with respect to an actualization since there is always a chance
--- i.e. except for the absurd property --- that it could be realized after the system has been changed
without destroying it.  It will become clear that a property can be actual or potential depending on the
state in which we consider the particular physical system.  Similarly we can say that certainty of obtaining
a positive answer when performing a definite experimental project depends on the state of the system.  In
order to construct our theory further we need to introduce the following relationship between definite
experimental projects, states and properties:
\bit  
\item
A definite experimental project $\alpha$ is called certain for a given singular realization of the given particular physical
system if it is sure that the positive response would be obtained should $\alpha$ be effectuated;
\item
A property $a$ is called actual for a given state if any, and so all, of the definite experimental projects corresponding to
the property $a$ are certain for any, and so all, of the singular realizations encoded by that state. A property is called
potential when it is not actual.
\eit
In
particular
is one of the essential achievements of sOQL that it gives a
consistent and coherent ontological account of physical properties contra
certain `overextrapolations' claimed to be inherent in quantum theory.
The
quantum mechanical formalism itself indeed allows a characterization of the properties of a quantum system as
being in
correspondence with the closed subspaces of the Hilbert space describing 
the
system in the above sense: Definite experimental projects $\alpha$ expressible in
quantum theory are of the form, ``the value of an observable is in
region
$E\subset\sigma(H)$'', where
$\sigma(H)$ is the spectrum of the self-adjoint operator $H$ describing the observable, thus we can write
$\alpha(H,E)$; more explicitly, the definite experimental project $\alpha(H,E)$ consists of measuring the
observable $H$ and obtaining an outcome in $E$;
the corresponding property $a$ is
then represented by the closed subspace of fixpoints of the projector $P_E$ that arises
via decomposition of $H$ according to von Neumann's spectral decomposition theorem since only
the states that imply the actuality of that property, that is, the states represented by a
ray
included in that subspace, will yield a positive outcome with certainty when ``we would perform
$\alpha(H,E)$''. 

\bigskip\par\noindent
Given the above notions, it becomes possible to introduce an operation on the collection ${\cal Q}$ of definite experimental
projects. What we have in mind is the product of a family of definite experimental projects $\Pi A$ which obtains its
operational meaning in the following way: 
\bit  
\item
The product  $\Pi A$ of a family $A$ of definite experimental projects is the definite experimental project
``choose arbitrarily one $\alpha$ in $A$ and effectuate it and attribute the obtained answer to $\Pi
A$\,''. 
\eit
More explicitly, given a particular realization of the system, $\Pi A$ is a certain definite experimental
project if and only if each member of $A$ is a certain definite experimental project.  Formally it
becomes possible to pre-order definite experimental projects by means of their certainty:
\bit 
\item
$\alpha \prec \beta$ := $\beta$ is certain whenever
$\alpha$ is certain.  
\eit
The notions of a trivial definite experimental project, which is always certain and an absurd
definite experimental project, which is never certain, can be introduced and play the role of respectively maximal and minimal
element of the collection
$Q$. The trivial and absurd definite experimental projects give rise on the level of properties to a trivial and absurd
property, the first is always actual while the latter is always potential.  Through the correspondence between definite
experimental projects and properties, the mentioned pre-order relation induces a partial order relation on
${\cal L}$\,:  
\bit 
\item
$a \leq b$ := $b$ is actual whenever $a$ is actual.    
\eit
The product of a collection of definite experimental projects with corresponding collection of properties
$A$ provides a greatest lower bound or meet for any $A\subseteq{\cal L}$\,, and as such the set of
properties
${\cal L}$ forms a complete lattice. Indeed, once we have the ``meet" of any collection of
properties, we can construct the operation of ``join" or least upper bound via Birkhoff's theorem stating
that for a given
$A
\subseteq {\cal L}$:
$$\bigvee A =
\bigwedge
\{x
\in {\cal L}\mid (\forall a \in A) a \leq x\}\,.$$  
Note however that contrary to the meet, the join admits of
no direct operational meaning in sOQL. Even more, we see that the join of a collection of properties need
not, and generically does not, correspond to a classical disjunction since the following implication is only secured
in one direction,
\bit 
\item
$(\exists a \in A$: $a$ is actual$)$ $\Rightarrow \bigvee A$ is actual.  
\eit
With respect to the meet we do have
\bit 
\item
$(\forall a\in A$: $a$ is actual$)$ $\Leftrightarrow \bigwedge A$ is actual,
\eit
following directly from the identification of the meet with the product of definite experimental projects.

\smallskip 
The property lattice description of a physical system in sOQL allows a
dual description by means of the system's states (Moore 1995, 1999) in terms of maximal state sets
$\mu(a)\subseteq
\Sigma$ for which a common property $a\in{\cal L}$ is actual whenever the system is in a state in
${\cal E}\in \mu(a)$\,, these sets being ordered by inclusion.  More explicitly, the state-property
duality may be straightforwardly characterized once we introduce the forcing relation
$\triangleright$ defined by ${\cal E} \triangleright a$ if and only if the property $a$ is actual in the state ${\cal E}$,
by   the
fact that we can associate to each property 
$a$ the set
$\mu(a) =
\{{\cal E} \in \Sigma \mid {\cal E} \triangleright a\}$ of states in which it is actual, and to each state ${\cal E}$ the set
$S({\cal E}) =
\{ a
\in {\cal L} \mid {\cal E} \triangleright a\}$ of its actual properties.  Formally $\mu: {\cal L} \to P(\Sigma): a \mapsto
\mu(a)$ is injective, satisfies the condition $\mu(\bigwedge A) = \bigcap \mu[A]$ and is called the Cartan map. In particular
we now see that actuality can be studied at the level of either states or properties since we have
$${\cal E} \in
\mu(a)
\Leftrightarrow {\cal E} 
\triangleright a \Leftrightarrow a \in S({\cal E})\,.$$  We will be a bit more explicit about the dual
description of a physical system by means of its states and operationally motivate the introduction of a
symmetric and antireflexive orthogonality relation on the set
$\Sigma$: 
\bit 
\item
Two states ${\cal E}$ and ${\cal E}'$ are called orthogonal, written ${\cal E} \perp {\cal E}'$, if there exists a definite
experimental project which is certain for the first and impossible for the second.
\eit   
If we equip $\Sigma$ with $\perp$, we can consider ${\cal
L}_{\Sigma}$ as the set of biorthogonal subsets, i.e., those $A \subseteq \Sigma$ with 
$A^{\perp \perp} = A$ for
$A^{\perp} = \{ {\cal E}'\mid(\forall {\cal E} \in A) {\cal E} \perp {\cal E}'\}$\,. 
Whenever $\perp$ 
is {\it separating}\, (Coecke et al (nd)a),  the codomain restriction of the Cartan map $\mu$ to the set of
biorthogonals gives us an isomorphism of complete atomistic orthocomplemented lattices (see below) extending
the dual description of a physical system by its state space 
$\Sigma$ and its property lattice
${\cal L}$ --- for details we refer to Moore (1999) and Coecke et al ((nd)a). It thus follows that ${\cal 
L}_{\Sigma}$ is an atomistic lattice where the singletons
$\{{\cal E}\}=\{{\cal E}\}^{\perp\perp}$ are the atoms.  Let us first introduce the notion of an atomistic
lattice more explicitly: 
\bit  
\item   
A complete lattice ${\cal L}$ is called atomistic if each element $a \in {\cal L}$ is generated by its 
subordinate atoms,
$$a = \bigvee\{p \in \Sigma_{\cal L} \mid p \leq a\}\,,$$ where the $p \in \Sigma_{\cal L}$ are by definition
the minimal nonzero elements of ${\cal L}$\,.
\eit
Under the assumption that states are in bijective correspondence with atoms, 
$\Sigma = \Sigma_{\cal L} \subseteq {\cal L}$\,,
each property lattice is atomistic in the sense that 
$$a = \bigvee\{ \bigwedge S({\cal E}) \mid {\cal E} \in \mu(a)\}$$
for each
$a \in {\cal L}$\,. Note that $\bigwedge S({\cal E})$ is the strongest actual property of the collection
$S({\cal E})$ which as such ``represents" the state.  

\smallskip
While we will only work with complete atomistic lattices in the following, we want to finish off this section with an
example explaining that a property lattice description for a quantum system will not lead to a boolean algebra. 
Under the assumption that properties have opposite properties and even more that each property $a \in
{\cal L}$ is the opposite of another one, we can formally introduce an orthocomplementation, i.e.:
\bit
\item 
A
surjective antitone involution
$': {\cal L} \to {\cal L}$ satisfying $a \wedge a' = 0$\,, ${a \leq b} \Rightarrow b' \leq a'$ and
$a''=a$\,. A lattice equipped with an orthocomplementation is usually abbreviated as an
ortholattice.
\eit
Consider now the property
lattice description of a photon as presented in Piron (1978)\,.  Take as a physical system a propagating photon which
is linearly polarized.  A definite experimental project $\alpha_\phi$ is
then defined by:
\bit
\item 
i. The apparatus: A polarizer oriented with angle $\phi$ and a counter placed behind it;
ii. The manual: Place the polarizer and counter within the passage of the photon;
iii. If one registers the passage of the photon through the polarizer, the result is ``yes"
and otherwise ``no".
\eit
\par\noindent {\hbox{\vrule height 3.7 truecm width 0pt
\special{illustration 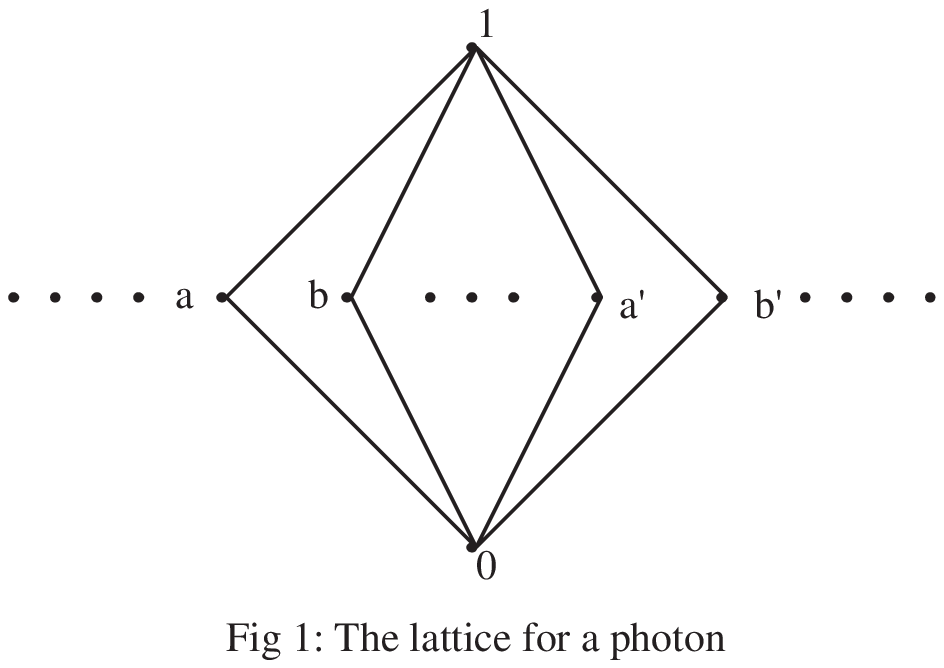 scaled 500}
\vrule height 0pt width 1 truecm}} \par
\par\vglue -3.8 truecm
\hangindent= 5.5 truecm \hangafter=-8
\noindent
Clearly, a property $a$ corresponding to $\alpha_\phi$ is called actual if it is
certain that $\alpha_\phi$ would lead to the response ``yes", should we perform the
experiment. In the diagram (Fig 1) of this photon, we consider some of its properties
explicitly: $a$ corresponds to
$\alpha_\phi$, $b$ to $\alpha_\phi'$, $a'$ to $\alpha_\phi + \pi/2$ and $b'$ to $\alpha_\phi' +
\pi/2$.  One immediately sees that this lattice is not boolean since distributivity is violated in the sense
that
$a\wedge (b\vee a') = a$ while $(a \wedge b) \vee (a \wedge a') = 0$. 

\vfill\eject\noindent  
{\bf 3. DYNAMIC OPERATIONAL QUANTUM LOGIC}  

\medskip\noindent 
Contrary to the static approach outlined above we  will now analyze how an actual property before
alteration will induce a property to be actual afterwards and, conversely, we characterize causes for
actuality.  The obtained result will then give rise to DoQL, when passing via
IoQL. We stress that both these approaches are still under full development. 
These developments were preceded by a representation theorem for deterministic evolutions of quantum 
systems as
given in Faure et al (1995) and for which DoQL provides an extension to non-deterministic cases. 
The new primitive concept (as compared to sOQL) is the notion of induction, defined in Amira et al (1998)
--- see also Coecke et al (2001) and Smets (2001)\,:  
\bit 
\item 
An induction $e \in \epsilon_s$ is a physical procedure that can be effectuated on a particular physical system $s$.  This
procedure, when carried out, might change $s$, modify the collection of its actual properties and thus its
state, or even destroy 
$s$.
\eit
On the collection $\epsilon_s$ of  
all inductions performable on a physical system 
$s$ we can consider two operations, one corresponding to the arbitrary choice of inductions and one corresponding to a finite
concatenation of inductions. Following
Amira et al (1998) and Coecke et al (2001)\,, the finite concatenation of inductions $e_1,e_2,..., e_n \in
\epsilon_s$ is the induction $e_1 \& e_2
\& ... \& e_n$ which consists of first performing $e_1$ then $e_2$, then ... until $e_n$. The arbitrary choice of
 inductions
in 
$\{e_i \M i \in I\} \subseteq \epsilon_s$ is the induction $\bigvee_i e_i$ consisting of performing one of
the $e_i$, chosen in any possible way.$^6$  Defined as such, an act of induction can for example be the assurance of a
free evolution (Faure et al 1995), indeterministic evolution (Coecke and Stubbe 1999), e.g., a measurement (Coecke and
Smets 2000, 2001), or the action of one subsystem
in a compound system on another one (Coecke 2000)\,. For reasons of formal simplicity we will only focus
on inductions which cannot lead to the destruction of the physical system under
consideration.$^7$  As such we presume that an induction cannot alter the
nature of a physical system, whereby we mean that what can happen is that a system's state is shifted
within the given initial state space, i.e., that the actuality and potentiality of the properties in the
initial property lattice is changed. This implies that the description of a physical system by means of its
state space or property lattice encompasses those states in which the system may be after performing an
induction. More explicitly, we can point
to the particular type of properties which are actual or potential before the system is altered and which contain
information about the actuality or potentiality of particular properties after the system is altered.  We
introduce this particular type of properties formally (Coecke et al 2001)\,:$^8$
\begin{equation}\label{inductionaction}
(e,a):\epsilon\times{\cal L}^{op}\to{\cal L}^{op}:(e,a)\mapsto e.a    
\end{equation}
where the reason for reversal of the lattice order (this is what ``$^{op}$'' in ${\cal L}^{op}$ stands
for) will be discussed below. The property $e.a$ stands for ``guaranteeing the actuality of a".  Existence
of a property $e.a$ is indeed operationally assured in the following way: 
given that $a$ corresponds to $\alpha$\,, $e.a$ corresponds to a definite experimental project
``$e.\alpha$'' of the form ``first execute the induction $e$ and then perform the definite experimental
project 
$\alpha$\,, and, attribute the outcome of $\alpha$ to $e.\alpha$"
(Faure et al 1995, Coecke 2000, Coecke et al 2001)\,.  In terms of actuality, following Smets (2001)\,:
\begin{quote}
``$e.a$ is an actual property for a system in a certain realization if it is sure that $a$
would be an actual property of the system should we perform induction $e$\,.''
\end{quote}  
This explains that if $e.a$ is actual it indeed ``guarantees" the actuality of $a$ with 
respect to $e$, while if $e.a$ is
potential it does not.  This expression crystallizes into the idea of introducing a 
{\it causal relation}:
$$\stackrel{e}{\leadsto}\ \subseteq {\cal L}_1 \times {\cal L}_2\,,$$  
where subscript $1$ points to the
lattice before $e$ and
$2$ to the lattice after $e$\,, as follows (Coecke et al 2001)\,:
\bit    
\item  
$a \stackrel{e}{\leadsto} b$ := the actuality of $a$ before $e$ induces (or, guarantees) the actuality of $b$ after
$e$.
\eit
Against the background of our characterization of $e.a$ we 
now see that $e.a \stackrel{e}{\leadsto} a$ will always be valid
and that
$$a \stackrel{e}{\leadsto} b\ \ \Longleftrightarrow\ \ a \leq e.b\,,$$  
so $\stackrel{e}{\leadsto}$ fully characterizes the action of $e.-:{\cal L}_2\to{\cal L}_1$\,. 
In case $e$ stands for the induction
``freeze" (with obvious significance, given a referential), conceived as timeless, then
$\stackrel{e}{\leadsto}$ reduces  to the partial
ordering $\leq$ of ${\cal L}_1 = {\cal L}_2$\,.       

\smallskip
To link the physical-operational level to a mathematical level, we associate with every 
induction $e \in \epsilon_s$ a map
called property propagation and a map called property causation (Coecke et al 2001)\,:
\smallskip\par\noindent
1) Property Causation:
\beqa
\bar{e}_\ast: {\cal L}_2 \to {\cal L}_1 : a_2 \mapsto e.a_2= \bigvee\{a_1 \in {\cal L}_1\mid a_1
\stackrel{e}{\leadsto} a_2 \}\,;
\eeqa
2) Property Propagation:
\beqa
\bar{e}^*:{\cal L}_1 \to {\cal L}_2: a_1 \mapsto \bigwedge\{a_2 \in {\cal L}_2 \mid a_1
\stackrel{e}{\leadsto} a_2\}\,.
\eeqa
Given those mappings, clearly $\bar{e}_{\ast}(a_2)$ is the weakest property whose actuality guarantees the
actuality of $a_2$ and
$\bar{e}^*(a_1)$ is the strongest property whose actuality is induced by that of $a_1$.  Further we
immediately obtain the Galois adjunction$^9$
$\bar{e}^*(\bar{e}_\ast(a_2)) \leq a_2$ and $a_1 \leq \bar{e}_{\ast}(\bar{e}^*(a_1))$\,, denoted as
$\bar{e}^*\dashv
\bar{e}_{\ast}$\,, since
\beqa  
a \leq \bar{e}_{\ast}(b)\ \ \Longleftrightarrow\ \ a \stackrel{e}{\leadsto} b\ \ \Longleftrightarrow\
\ \bar{e}^*(a)
\leq b\,.
\eeqa
In Coecke et al (2001) this adjunction is referred to as ``causal duality'', since it expresses the dual
expressibility of dynamic behavior for physical systems respectively in terms of propagation of properties
and causal assignment.$^{10}$ We also  recall here that this argument 
towards causal duality suffices to establish evolution for quantum systems, i.e., systems with the lattice
of closed subspaces of a Hilbert space as property lattice, in terms of linear or anti-linear maps (Faure et
al 1995) and compoundness in terms of the tensor product of the corresponding Hilbert spaces (Coecke
2000)\,. It also follows from the above that the action defined in eq.(\ref{inductionaction}) defines a
quantale module action (Coecke et al 2001) --- quantales will be discussed below. Crucial here is the fact that
$(\bigvee_ie_i)\cdot\alpha$ and $\Pi_i(e_i\cdot\alpha)$ clearly define the same property
$(\bigvee_ie_i)\cdot a=\bigwedge_i(e_i\cdot a)$\,, since both $\Pi$ for definite experimental projects and
$\bigvee$ for inductions express choice. Accordingly, the opposite ordering in eq.(\ref{inductionaction})
then matches ${\cal L}$-meets with $\epsilon$-joins.

\smallskip 
In the above we discussed the strongest property $\bar{e}^*(a_1)$ of which actuality is induced by an
induction due to actuality of $a_1$\,, but only for maximally deterministic evolutions this fully
describes the system's behavior. In other cases it makes sense to consider how (true logical) disjunctions
of properties propagate, as such allowing accurate representation of for example the emergence of
disjunction in a perfect quantum measurement (Coecke 2002, Coecke and Smets 2000, Smets
2001) due to the uncertainty on the measurement outcome whenever the system is not in an
eigenstate of this measurement.$^{11}$  First we introduce
the notion of an actuality set as a set of properties of which at least one element is actual, clearly
encoding logical disjunction in terms of actuality.  As before we want to express propagation and
causal assignment of these actuality sets.  While we shift from the level of properties to sets of
properties, we will have to take care that we don't loose particular information on the structure of
${\cal L}$\,, in particular its operationally derived order.  
The solution to this problem consists in considering a certain kind of ideal. More specifically we work
with the set
$DI({\cal L})
\subseteq P({\cal L})$\,, $P({\cal L})$ being the powerset of ${\cal L}$\,, of which the elements are called
property sets and which formally are the so-called distributive ideals of ${\cal L}$\,, introduced
in a purely mathematical setting in Bruns and Lakser (1970)\,:
\bit 
\item
A distributive ideal is an order ideal, i.e. 
if $a \leq b \in I$ then $a \in I$ and $I \not = \emptyset$, and is closed
under distributive joins, i.e. if $A \subseteq I \in DI({\cal L})$ then $\bigvee A \in I$ whenever 
we have $$\forall b \in
{\cal L}: b
\wedge
\bigvee A = \bigvee\{b \wedge a \mid a \in A \}.$$  
\eit
Intuitively, this choice can be motivated as follows (a much more rigid argumentation does
exists):$^{12}$
i. a first choice for encoding disjunctions would be the powerset itself, however, if $a\leq b$ we don't
have $\{a\}\subseteq\{b\}$ so we do not preserve order; otherwise stated, if $a<b$ then the ``propositions''
$\{a\}$ and $\{a,b\}$ (read: either $a$ \underline{or} $b$ is actual) mean the same thing, since actuality
of $b$ is implied by that of $a$\,; ii. we can clearly overcome this problem by considering order ideals
$$I({\cal L}):=\{\downarrow\![A]|A\subseteq{\cal L}\}\subset P({\cal L})\,;$$ however, in
case the property lattice would be a complete Heyting algebra in which
all joins encode disjunctions, then $A$ and $\{\bigvee A\}$ again mean
the same thing; this redundancy is then exactly eliminated by considering
distributive ideals (Coecke 2002, Coecke and Smets 2001)\,.  For ${\cal L}$ atomistic and
$\Sigma \subseteq {\cal L}$\,, $DI({\cal L}) \cong P(\Sigma)$ which
implies that $DI({\cal L})$ is a complete atomistic boolean algebra
(Coecke 2002)\,. 

\smallskip 
Similarly as for properties, for property sets we can operationally motivate the introduction of a causal
relation
$\stackrel{e}{\leadsto}\ \subseteq\, DI({\cal L})_1 \times DI({\cal L})_2$\,:
\bit 
\item
$A \stackrel{e}{\leadsto} B$ := if property set $A$ is an actuality set before $e$\,, $A$ induces that
property set
$B$ is an actuality set after $e$\,.
\eit
To every induction we associate a map called
property set causation and a map called property set propagation:
\smallskip\par\noindent
1) Property set causation
\beqa
\hat{e}_{\ast} : DI({\cal L})_2 \to DI({\cal L})_1: A_2 \mapsto {\cal
C}(\bigcup\{A_1 \in DI({\cal L})_1\mid A_1 \stackrel{e}{\leadsto} A_2\})
\eeqa 
2) Property set propagation
\beqa
\hat{e}^* : DI({\cal L})_1 \to DI({\cal L})_2: A_1 \mapsto \bigcap \{A_2 \in DI({\cal L})_2
\mid A_1
\stackrel{e}{\leadsto} A_2\}
\eeqa
where 
$${\cal C}: P({\cal L}) \to P({\cal L}): A \mapsto \bigcap\{B \in DI({\cal L})\mid A \subseteq B\}\,.$$
Similar as above we obtain an adjunction $\hat{e}^*\ \dashv\ \hat{e}_{\ast}$ following from:
\beqa  
\hat{e}^* (A)_1 \subseteq A_2\ \ \Longleftrightarrow\ \ A_1 \stackrel{e}{\leadsto} A_2\ \
\Longleftrightarrow\ \ A_1 \subseteq \hat{e}_{\ast}(A_2)
\eeqa 
In case the induction $e$ stands for ``freeze" we obtain $A \stackrel{e}{\leadsto} B = A \subseteq B$.
Note that the join preservation that follows from the adjunction $\hat{e}^*\dashv\hat{e}_{\ast}$ expresses
the physically obvious preservation of disjunction for temporal processes. 

\smallskip
It is however also important to stress here that not all maps $\hat{e}^* : DI({\cal L})_1 \to DI({\cal
L})_2$ are physically meaningful. Indeed, any physical induction admits mutually adjoint maps $\bar{e}^*$ and
$\bar{e}_*$ with the significance discussed above, the existence of such a map $\bar{e}^*:{\cal L}_1\to{\cal L}_2$\,
forcing $\hat{e}^*$ to satisfy a join continuity condition (Coecke and Stubbe 1999, Coecke et al 2001, Coecke 2002),
namely:
\beq\label{continuity}
\bigvee A=\bigvee B\ \ \Longrightarrow\ \ \bigvee\hat{e}^*(A)=\bigvee\hat{e}^*(B)\,,
\eeq  
which indeed expresses well-definedness of a corresponding $\bar{e}^*$ since given an actuality set
$A$\,, the strongest property that is actual with certainty is $\bigvee A$\,. It is exactly in the existence
of a non-trivial condition as in eq.(\ref{continuity}) that the non-classicality of quantum theory comes
in. As such, both causal dualities, the one on the level of properties and the other on the level of property sets, 
provide a {\it physical law} on transitions, respectively condition eq.(\ref{continuity}) and preservation of
$DI({\cal L})$-joins. We also want to stress a manifest difference here with the setting in van Benthem
(1994)\,:
\begin{quote}
``The most general model of dynamics is simply this: some system moves through a space of possibilities.
Thus there is to be some set [$\Sigma$] of relevant {\it states} (cognitive, physical, etc.) and a
family $[\{R_e|\,e\in\epsilon\}]$ of binary {\it transition relations} among them, corresponding to
actions that could be performed to change from one state to another. [...] Let us briefly consider a number
of dynamic `genres', [...] $\bullet$ Real action in the world changes actual physical states. [...] What
are most general operations on actions? Ubiquitous examples are {\it sequential composition} and {\it
choice}\,.'' (p.109, 110, 112)
\end{quote}
Thus it seems to us that the author aims to cover our study of dynamic behavior of physical systems. He
moreover states (van Benthem 1994)\,: 
\bq
``The main claim of this paper is that the above systems of relational algebra and dynamic logic provide
a convenient architecture for bringing out essential logical features of action and cognition.'' (p.130)
\eq  
We claim here that 
relational structures are inappropriate for modeling physical dynamics, {\it even classically}\,! Let us
motivate this claim. It follows from the above that transitions of properties of physical systems are
internally structured by the causal duality, which in the particular case of non-classical systems
restricts possible transitions.  Recalling that for atomistic property lattices we have $P(\Sigma)\cong
DI({\cal L})$\, it is definitely true that any join preserving map $\hat{e}^*:DI({\cal L})_1\to DI({\cal L})_2$
defines a unique relation  
$R_e\subseteq \Sigma_1\times\Sigma_2$\,, and conversely, any relation $R\subseteq
\Sigma_1\times\Sigma_2$ defines a unique (join preserving) map $f_R:DI({\cal L})_1\to DI({\cal L})_2$\,. Next, any
relation
$R\subseteq
\Sigma_1\times\Sigma_2$ has an inverse $R^{-1}\subseteq \Sigma_2\times\Sigma_1$ and this inverse plays a major role in van
Benthem (1991, 1994) as {\it converse action}\,.  However, nothing assures that when $\hat{e}^*$ satisfies
eq.(\ref{continuity}) that the map $f_{R_e^{-1}}:DI({\cal L})_2\to DI({\cal L})_1$ encoding the
relational inverse satisfies  eq.(\ref{continuity})\,, and as such, has any physical significance
at all.$^{13}$  Obviously, this argument applies only to non-classical systems.  More generally, however, since it is
the duality between causation and propagation at the
$DI({\cal L})$-level that guarantees preservation of disjunction we feel that it should be present in any
modelization, and although relations and union preserving maps between powersets are in bijective
correspondence, they have fundamentally different dual realizations: relations have inverses, and union
preserving maps between powersets have adjoints, and these two do not correspond at all, respectively being
encoded as (in terms of maps between powersets):
\beqa
f_{R_{e}^{-1}}(A)\!\!&=&\!\!{\{p\in\Sigma|\exists q\in A:q\in \hat{e}^*(\{p\})\}}\\
\hat{e}_*(A)\!\!&=&\!\!{\{p\in\Sigma|\forall q\in A:q\in \hat{e}^*(\{p\})\}}
\eeqa
if it was even only by the fact that one preserves unions and the other one intersections. As such, the
seemingly innocent choice of representation in terms of relations or union preserving maps between powersets
does have some manifest consequence in terms of the implementation of causal duality.

\smallskip
In the remaining part of this section we concentrate on the logic of actuality sets as initiated in 
Coecke (2002).  Introduce the following primitive connectives:
\beqa 
\bigwedge_{DI({\cal L})}& :& P(DI({\cal L})) \to DI({\cal L}): {\cal A} \mapsto \bigcap {\cal A}\,;\\
\bigvee_{DI({\cal L})}& :& P(DI({\cal L})) \to DI({\cal L}): {\cal A} \mapsto {\cal C}(\bigcup {\cal A})\,;\\
\rightarrow_{DI({\cal L})}& :& DI({\cal L}) \times DI({\cal L}): (B,C) \mapsto \bigvee_{DI({\cal L})}
\{A\in DI({\cal L})\mid A \cap B\subseteq C\}\\
&&\hspace{5.8cm}=\{a\in {\cal L}|\forall b\in B:a\wedge b\in C\}\,;\\ 
{\cal R}_{DI({\cal L})}& :& DI({\cal L}) \to DI({\cal L}): A \mapsto \downarrow( \bigvee_{\cal L} A)\,.
\eeqa
While the first three connectives are standard in intuitionistic logic, ${\cal R}_{DI({\cal L})}$ should be
conceived as a ``resolution-connective" allowing us to recuperate the logical structure of properties on
the level of property sets (Coecke 2002). In particular, for classical systems we have ${\cal R}_{DI({\cal
L})}=id_{DI({\cal L})}$.  Note that the condition in eq.(\ref{continuity}) now becomes:
\beq\label{continuity2}
{\cal R}_{DI({\cal L})}(A)={\cal R}_{DI({\cal L})}(B)\ \ \Longrightarrow\ \
{\cal R}_{DI({\cal L})}(\hat{e}^*(A))={\cal R}_{DI({\cal L})}(\hat{e}^*(B))\,, 
\eeq
restricting the physically admissible transitions. Clearly, this condition is trivially satisfied for
classical systems.

\smallskip  
When concentrating on the material implication, we want to stress that on the level of ${\cal L}$ only for
the properties in a distributive sublattice we can say the following: 
\beqa
p \models (a \rightarrow_{\cal L} b) &\ \Longleftrightarrow\ & \{p\} \cap \mu(a) \subseteq \mu(b)\\
&\ \Longleftrightarrow\ & p \in \mu(a) \Rightarrow p \in \mu(b)\\
&\ \Longleftrightarrow\ & p \models a \Rightarrow p \models b\,.
\eeqa
In that case, this implication satisfies the strengthened law of entailment: 
\beqa
\mu(a \rightarrow_{\cal L} b) = \Sigma\ \ \Longleftrightarrow\ \ \mu(a) 
\subseteq \mu(b)\ \ \Longleftrightarrow\ \ a \leq b\,. 
\eeqa
 Note that in the non-distributive case for $a \rightarrow_{\cal L} b = a' \vee b
\in {\cal L}$ we can only say that $(p \models a
\rightarrow_{\cal L} b) \Leftarrow (p \models a \Rightarrow p \models b)$, which goes together with the fact that there
are examples of orthomodular lattices for which $a' \vee b = 1$ while $a \not \leq b$.  Hence in general, for $T= \{p
\in
\Sigma \mid p \models a \Rightarrow p \models b\}$ there will not be an element $x \in {\cal L}$ for which $\mu(x) = T$.
There are of course examples of other implications which do satisfy the strengthened law of entailment in the orthomodular
case --- see for instance Kalmbach (1983)\,. It is now our aim to focus on the
implication for elements in $DI({\cal L})$.  First we lift the Cartan map to the level of property sets $\mu(A) :=
\bigcup
\mu[A] \subseteq \Sigma_{\cal L}$\,, then we obtain the following semantical interpretation:
\beqa
p \models (A \rightarrow_{DI({\cal L})} B) &\ \Longleftrightarrow\  & \{p\} \cap \mu(A) \subseteq \mu(B) \\
&\ \Longleftrightarrow\ & p \in \mu(A) \Rightarrow p \in \mu(B)\\
&\ \Longleftrightarrow\ & p \models A \Rightarrow p \models B\,.
\eeqa
From this it follows that 
$$\mu (A \rightarrow_{DI({\cal L})} B) = \{p \in \Sigma_{\cal L} \mid (p \models A) \Rightarrow (p
\models B)\}\,,$$ which allows us to reformulate the given static implication $(-\rightarrow_{DI({\cal L})}-)$
as follows (Coecke (nd), Coecke and Smets 2001):
\beqa 
(A\to_{{DI}({L})} B) \!\!&=&\!\! \bigvee_{{DI}({L})} \{ C \in {DI}({L}) \mid \forall D\vdash C:(D \vdash A
\Rightarrow D
\vdash B)\}\\ 
\!\!&=&\!\! \{ c\in L \mid \forall d\leq c:(d\in A \Rightarrow d\in  
B)\}\,,  
\eeqa
where $D \vdash A \Leftrightarrow \forall p \in \mu(D): p \in \mu(A)$.  Again this implication satisfies the
strengthened law of entailment in the sense that 
$$(A \rightarrow_{DI({\cal L})} B) = {\cal L}\ \ \Longleftrightarrow\ \ A \subseteq B$$
and is as
operation the right adjoint with respect to the conjunction $\bigwedge_{DI({\cal L})}$\,, as it is always
the case for the implication connective of a complete Heyting algebra --- see for example Borceux (1994)
or Johnstone (1982)\,:
$$(A\wedge_{DI({\cal L})}-)\dashv(A\rightarrow_{DI({\cal L})}-)\,.$$ 
  
\smallskip 
As for the static material implication defined above, we now want to look for a dynamic propagation-implication
satisfying the following:
\beqa
(A \stackrel{e}{\rightarrow} B) = {\cal L}\ \ \Longleftrightarrow\ \ A \stackrel{e}{\leadsto} B\,.
\eeqa
The candidate which naturally arises is (Coecke (nd), Coecke and Smets 2001):
\beqa 
(A \stackrel{e}{\rightarrow} B) : = \{ c\in L \mid \forall d\leq c:(d\in A \Rightarrow \hat{e}^*
(\downarrow\!d)\subseteq  B)\}\,.
\eeqa
indicated by a semantical interpretation as
\beqa 
\mu (A \rightarrow_{DI({\cal L})} B) = \{p \in \Sigma_{\cal L} \mid (p \models A) \Rightarrow (\hat{e}^* (\{p\})
\models B)\}\,.
\eeqa
In case $e$ stands for the induction ``freeze" we see that $\stackrel{e}{\rightarrow}$ reduces to 
$\rightarrow_{DI({\cal L})}$\,.  Similar as in the static case, we can find an induction-labeled operation
as left adjoint to the dynamic propagation-implication, i.e., $$(A\otimes_e-)\dashv
(A\stackrel{e}{\rightarrow}-)\ \ \
\mbox{\rm with}\ \ \ A \otimes_e B := \hat{e}^* (A \wedge_{DI({\cal L})} B)\,.$$   
It is important here to
notice that this dynamic conjunction is a commutative operation.  Since $\hat{e}^*$ preserves joins and
since in $DI({\cal L})$ binary meets distribute over arbitrary joins (being a complete Heyting algebra) we
moreover have (Coecke (nd))\,: 
$$A\otimes_e (\bigvee_{DI({\cal L})} {\cal B}) =
\bigvee_{DI({\cal L})} \{A \otimes_e B \mid B \in {\cal B}\}\,,$$ 
so $DI({\cal L})$ is equipped with operations $\stackrel{e}{\rightarrow}$ and $\otimes_e$\,, for every $e
\in 
\epsilon_s$ the latter yielding ``commutative quantales''.  Let us give the definition of a
quantale (Rosenthal 1990, 1996, Paseka and Rosicky 2000)\,:  
\bit 
\item 
A {\it quantale} is a complete lattice ${Q}$ together with an associative binary 
operation $\circ$ that satisfies 
$a \circ (\bigvee_i b_i)= \bigvee_i(a \circ b_i)$ and $(\bigvee_i b_i) \circ a = \bigvee_i (b_i \circ
a)$ for all $a, b_i\in{Q}$\,.
\eit
Thus, for each induction $e$ we obtain that $(DI({\cal L}), \bigvee_{DI({\cal L})}, \otimes_e)$ is a
commutative quantale since $\otimes_e$ is a commutative operation.  In case $e$ stands for ``freeze", the
mentioned quantale becomes a locale
$(DI({\cal L}), \bigvee_{DI({\cal L})},\wedge_{DI({\cal L})})$\,, i.e., a complete Heyting algebra. Recall
here that a locale is a quantale with as quantale product the meet-operation of the complete lattice, and
one verifies that this definition exactly coincides with that of a complete Heyting algebra --- for details
we refer again to Borceux (1994) or Johnstone (1982)\,.

\smallskip  
As for the propagation-implications which, when valid, express a forward causal relation between property sets, we
can introduce causation-implications.  The relation to which these causation-implications match will be a backward
relation introduced as follows (Coecke and Smets 2001):
\bit 
\item
$A \stackrel{e}{\looparrowleft} B$ := If property set $B$ is necessarily an actuality set after $e$ then property set
$A$ was an actuality set before $e$\,.
\eit
Formally we see  
$$A \stackrel{e}{\looparrowleft}B\ \ \Longleftrightarrow\ \ \hat{e}_\ast(B) \subseteq A\,.$$  
The causation-implications we want to work with now have to satisfy:
\beqa 
A \stackrel{e}{\leftarrow} B = {\cal L}\ \ \Longleftrightarrow\ \ A \stackrel{e}{\looparrowleft} B\,.
\eeqa  
The candidate which satisfies this condition is:
\beqa 
(A \stackrel{e}{\leftarrow} B) : = \{ c\in L \mid \forall d\leq c:(d\in A \Leftarrow \hat{e}^*
(\downarrow\!d)\subseteq   B)\}\,.
\eeqa
indicated by a semantical interpretation as
\beqa 
\mu (A \leftarrow_{DI({\cal L})} B) = \{p \in \Sigma_{\cal L} \mid (p \models A) \Leftarrow (\hat{e}^* (\{p\})
\models B)\}\,.
\eeqa
As such we see that when $A \stackrel{e}{\leftarrow}B$ is valid (i.e. equal to ${\cal L}$) then it
expresses that if property set $B$ is an actuality set after $e$ then property set $A$ was an actuality set
before $e$. Again we have a left adjoint for the causation-implication: $$(-
{\,_e\!\otimes} B) \dashv (- \stackrel{e}{\leftarrow}B)\ \ \
\mbox{\rm with}\ \ \ A {\,_e\!\otimes} B := A \wedge_{DI({\cal L})}  \hat{e}_\ast (B)\,.$$ 
Thus we can additionally equip
$DI({\cal L})$ with $\stackrel{e}{\leftarrow}$ and ${_e\otimes}$\,, for each $e \in \epsilon_s$ the latter
yielding non-commutative co-quantales $(DI({\cal L})\,,
\bigwedge_{DI({\cal L})}\,, {_e\otimes})$\,, i.e., with a distributive property with respect to meets.  
Note that the preservation of joins for propagation versus that of meets for causation reflects
here in a two-sided distributivity respectively with respect to joins and meets. Indeed, since we have:
$$
({\cal L}\otimes_e -) = \hat{e}^* (-)\quad\quad\quad\quad  
({\cal L}{\,_e\!\otimes} -) = \hat{e}_\ast (-)
$$
this distributivity truly encodes the respective join and meet preservation and consequently,
also the causal duality.  

\smallskip
It is important to remark that the semantics we obtain is the complete Heyting algebra of
actuality sets 
$DI({\cal L})$ equipped with additional dynamic connectives to express causation and propagation:
\beqa
\Bigr(DI({\cal L})\,, {\cal R}_{DI({\cal L})}\,,\bigvee_{DI({\cal L})}\,,
\bigwedge_{DI({\cal L})}\,,\neg,\,\{ \otimes_e,{_e \otimes},\stackrel{e}{\to},\stackrel{e}{\leftarrow},
\stackrel{e}{\neg}\M e \in
\epsilon \}\Bigl)\,.
\eeqa
It turns out that we have a forward negation $\neg$ which does not depend on $e$\,, and thus
coincides with that of ``freeze'', i.e., the static one, and a backward negation $\stackrel{e}{\neg}$ which by
contrast does depend on $e$ --- for a discussion of these negations we refer to Coecke (nd) and Smets (2001)\,. 
Note here that $DI({\cal L})$ is also a left quantale module for 
$\hat{\epsilon} =
\{\hat{e}\,|\,e \in
\epsilon_s\}$ when considering the pointwise action of $(e.-)$\,. This then intertwines the two quantale
structures that emerge in our setting. We end this paragraph by mentioning that it is possible to
implement other kinds of implications on
$DI({\cal L})$ that extend the causal relation. As an example, it is
possible to extend 
$DI({\cal L})$ with bi-labeled families of non-commutative quantales rendering bi-labeled implications,
some of them extending the ones presented in this paper.

\smallskip
In Coecke and Smets (2001) it is argued that the Sasaki adjunction is an incarnation of causal duality for the
particular case of a quantum measurement $\varphi_a$ with a projector (on $a$) as corresponding self-adjoint
operator. 
This has as a striking consequence, since validity of the Sasaki adjunction is equivalent
to ``orthomodularity", sOQL embodies a hidden dynamical ingredient which is
algebraically identifiable as orthomodularity.  One could as such argue that the necessity of the passage from sOQL
to DoQL was already announced within sOQL itself, it was just waiting to be revealed.  As a more radical statement
one could say that due to this hidden dynamical ingredient, it is impossible to give a full sense to
quantum theory in logical terms within an essentially static setting.  Following Coecke and Smets (2001),
this fact can be deduced from DoQL by eliminating the emergent disjunctivity when introducing modalities
with respect to actuality and conditioning. We can in that case derive that the labeled dynamic hooks
that encode quantum measurements act on properties as
\beqa
\ \ \ \ (a_1\stackrel{\varphi_a}{\to}a_2):= 
(a_1\to_{\cal L}(a\stackrel{S}{\to} a_2))\quad {\rm and}
\quad (a_1\stackrel{\varphi_a}{\leftarrow}a_2):=
((a\stackrel{S}{\to}a_2)\to_{\cal L} a_1)
\eeqa 
where $(-\stackrel{S}{\to}-)$ is the well-known Sasaki hook and we identify $a$ and $\{a\}$. One could say
that the transition from  either classical or  intuitionistic logicality to 
``true''quantum logicality entails besides the introduction of an
additional unary connective ``operational resolution'' the shift from a binary implication connective to a
ternary connective where two of the arguments have an ontological connotation and the third,  
the new one, an empirical.

\bigskip\noindent   
{\bf 4. COMPARISON WITH LINEAR LOGIC}

\medskip\noindent
We will analyze another logic of dynamics, namely linear logic, while focusing on the differences with DoQL
especially with respect to the quantale structures mentioned above. We intend to give a brief overview of the basic
ideas behind ``Linear logic" as introduced in Girard (1987, 1989)\,.  It's categorical semantics in terms of
a $^*$-autonomous category appeared in Barr (1979) and it is fair to say that already in Lambek (1958) a
non-commutative fragment of linear logic was present. We follow the discussion of Smets (2001).  

\smallskip
The main advantage  linear logic has with respect to classical/intuitionistic logic is that
it allows us to deal with actions versus situations in the sense of stable truths (Girard 1989).  This should be understood
in the sense that linear logic is often called {\it a resource sensitive logic}.  The linear logical formulas can be
conceived as expressing finite resources, the classical formulas then being interpretable as corresponding to
unlimited or eternal resources.   Allowing
ourselves to be a bit more formal on this matter, resource sensitivity is linked to the explicit control of the
weakening and contraction rules.  As structural rules, weakening and contraction will be discarded in the general
linear logical framework.  Note that in non-commutative linear logic, to which we will come back later, the
exchange-rule will also be dropped. We use
$\mathfrak{A},\mathfrak{B}$ for sequences of well formed formulas and
$\mathfrak{a},\mathfrak{b}$ for well formed formulas. Sequents are conceived as usual: 
\beqa  
\begin{array}{ccccc}
(weakening L) & {\mathfrak{A} \pijl \mathfrak{B}\over \mathfrak{A},\mathfrak{a} \pijl \mathfrak{B}} & &
(weakening R) & {\mathfrak{A} \pijl \mathfrak{B} \over \mathfrak{A} \pijl \mathfrak{B},\mathfrak{a}}  \\ \\
(contraction L) & {\mathfrak{A}, \mathfrak{a}, \mathfrak{a} \pijl \mathfrak{B} \over \mathfrak{A},\mathfrak{a} \pijl \mathfrak{B}}&& 
(contraction R) & {\mathfrak{A} \pijl \mathfrak{a},\mathfrak{a}, \mathfrak{B} \over \mathfrak{A} \pijl \mathfrak{a}, \mathfrak{B}}\\ \\
(exchange L) & {\mathfrak{A}_1, \mathfrak{a},\mathfrak{b},\mathfrak{A}_2 \pijl \mathfrak{B} \over \mathfrak{A}_1,\mathfrak{b},\mathfrak{a}, \mathfrak{A}_2 \pijl
\mathfrak{B}} &&  (exchange R) & {\mathfrak{A} \pijl \mathfrak{B}_1,\mathfrak{a},\mathfrak{b},\mathfrak{B}_2 \over \mathfrak{A} \pijl
\mathfrak{B}_1,\mathfrak{b},\mathfrak{a},\mathfrak{B}_2}
\end{array}
\eeqa
Dropping weakening and contraction implies that linear formulas cannot be duplicated or contracted at
random, in other words, our resources are restricted. An important consequence of dropping these two
structural rules is the existence of two kinds of ``disjunctions'' and ``conjunctions''.$^{14}$ We will obtain a
so-called additive disjunction
$\oplus$ and additive conjunction
$\sqcap$ and a so-called multiplicative disjunction $\acht$ and multiplicative conjunction $\otimes$.  The following left
and right rules will make their differences clear.
\beqa
\begin{array}{lllll}

(\sqcap\ R) & {\mathfrak{A} \rightarrow \mathfrak{a}, \mathfrak{B} \quad \mathfrak{A} \pijl \mathfrak{b}, \mathfrak{B} \over
\mathfrak{A} \pijl \mathfrak{a} \sqcap
\mathfrak{b}, \mathfrak{B}} & &

(\oplus\ L) & {\mathfrak{A}, \mathfrak{a} \pijl \mathfrak{B} \quad \mathfrak{A}, \mathfrak{b} \pijl \mathfrak{B} \over \mathfrak{A}, \mathfrak{a} \oplus
\mathfrak{b} \pijl \mathfrak{B}}\\ \\

(\sqcap\ L) & {\mathfrak{A}, \mathfrak{a} \pijl \mathfrak{B} \over \mathfrak{A}, \mathfrak{a} \sqcap \mathfrak{b} \pijl
\mathfrak{B}} \quad {\mathfrak{A},\mathfrak{b}
\pijl \mathfrak{B} \over \mathfrak{A}, \mathfrak{a}\sqcap \mathfrak{b} \pijl \mathfrak{B}} & &

(\oplus\ R) & {\mathfrak{A} \pijl \mathfrak{a}, \mathfrak{B} \over \mathfrak{A} \pijl \mathfrak{a} \oplus \mathfrak{b}, \mathfrak{B}} \quad {\mathfrak{A}
\pijl \mathfrak{b}, \mathfrak{B} \over \mathfrak{A} \pijl \mathfrak{a} \oplus \mathfrak{b}, \mathfrak{B}}\\ \\

(\otimes\ R) & {\mathfrak{A}_1 \pijl \mathfrak{a},\mathfrak{B}_1 \quad \mathfrak{A}_2 \pijl \mathfrak{b}, \mathfrak{B}_2 \over \mathfrak{A}_1,\mathfrak{A}_2
\pijl
\mathfrak{a}
\otimes
\mathfrak{b},\mathfrak{B}_1,\mathfrak{B}_2} & & 

(\acht\ L) & {\mathfrak{A}_1, \mathfrak{a} \pijl \mathfrak{B}_1 \quad \mathfrak{A}_2, \mathfrak{b} \pijl \mathfrak{B}_2 \over \mathfrak{A}_1,\mathfrak{A}_2,
\mathfrak{a}\acht \mathfrak{b} \pijl \mathfrak{B}_1,\mathfrak{B}_2} \\ \\

(\otimes\ L) & {\mathfrak{A},\mathfrak{a},\mathfrak{b} \pijl \mathfrak{B} \over \mathfrak{A}, \mathfrak{a} \otimes \mathfrak{b} \pijl \mathfrak{B}} & &

(\acht\ R) & {\mathfrak{A} \pijl \mathfrak{a}, \mathfrak{b}, \mathfrak{B} \over \mathfrak{A} \pijl \mathfrak{a} \acht \mathfrak{b},\mathfrak{B}}

\end{array}
\eeqa
By allowing the structural rules and by using $(\otimes)$ we can express the
$(\sqcap)$-rules and vice versa. A similar result can be obtained for the $(\oplus)$ and $(\acht)$-rules.

\smallskip
To understand linear logic, it is necessary to take a look at the intuitive meaning of the above
additives and multiplicatives derived from their use by the above rules.  First it is important
to note that to obtain
$\otimes$ in a conclusion, no sharing of resources is allowed, while the contrary is the case for $\sqcap$. Similarly, there
is a difference between $\acht$ and $\oplus$.  In the line of thought exposed in Girard (1989), the meanings
to be attached to the connectives are the following:
\bit
\item 
$\mathfrak{a} \otimes \mathfrak{b}$ means that both resources, $\mathfrak{a}$ and $\mathfrak{b}$ are
given simultaneously\,; 
\item 
$\mathfrak{a} \sqcap \mathfrak{b}$ means that one may choose between $\mathfrak{a}$ and $\mathfrak{b}$\,;
\item 
$\mathfrak{a} \oplus \mathfrak{b}$ means that one of both resources, $\mathfrak{a}$  or $\mathfrak{b}$, is given though we
have lack of knowledge concerning the exact one\,;
\item
$\mathfrak{a} \acht \mathfrak{b}$ expresses a constructive disjunction.
\eit
The meaning of  $\acht$ becomes clearer {\it
when} we follow J.Y. Girard in his construction that every atomic formula of his linear logical language has by definition
a negation $(-)\lood$.  Running a bit ahead
of our story, the meaning of $\mathfrak{a} \acht \mathfrak{b}$ will now come down to the situation where ``if" not
$\mathfrak{a}$ is given ``then" $\mathfrak{b}$ is given and ``if" not $\mathfrak{b}$ is given ``then" $\mathfrak{a}$ is
given.  Of course this explanation is linked to the commutative case where a linear logical implication is {\it defined} as 
$\mathfrak{a}\lood \acht \mathfrak{b} :=
\mathfrak{a}$ $\multimap$ $\mathfrak{b}$ which by transposition equals $\mathfrak{b}\lood$ $\multimap$ $\mathfrak{a}\lood$. 
In the same sense, only in the commutative case where $\mathfrak{a} \otimes \mathfrak{b}$ equals $\mathfrak{b} \otimes
\mathfrak{a}$, does it make sense to say that
$\mathfrak{a} \otimes \mathfrak{b}$ comes down to {\it simultaneous} given resources.

\smallskip
We will be more specific on the linear logical implication and
analyze the underlying philosophical ideas as presented in Girard (1989) where of course $\multimap$ is
defined by means of
$\acht$. What is important is that the linear implication should mimic exactly what happens when a non-iteratable action
is being performed, where we conceive of a non-iteratable action to be such that after its performance the initial
resources are not available any more {\it as} initial resources.  The linear implication should as such express the
{\it consumption of initial resources} and simultaneously the {\it production of final resources}. Indeed, as stated
in Girard (1995), the linear implication $\multimap$ expresses a form of causality:  
$\mathfrak{a}$ $\multimap$ $\mathfrak{b}$ is to be conceived as ``from $\mathfrak{a}$ get $\mathfrak{b}$".  
More explicitly (Girard 1989)\,:
\bq 
``A causal implication cannot be iterated since the conditions are modified after its use; this process of modification of
the premises (conditions) is known in physics as {\it reaction}." (p.72)
\eq 
or in Girard (2000):
\bq
``C'est donc une vision {\it causale} de la d\'eduction logique, qui s'oppose \`a la p\'erennit\'e de
la v\'erit\'e traditionelle en philosophie et en math\'e-matiques. On a ici des v\'erit\'es fugaces,
contingentes, domin\'ees par l'id\'ee de ressource et d'action.'' (p.532)
\eq
If we understand this correctly, the act of
consumption and production is called a non-iteratable action while the process of modification
of initial conditions, the deprivation of resources, is called reaction.  The idea of relating action and reaction is in a
sense metaphorically based on Newton's action-reaction principle in physics.  Girard uses this metaphor also when he
explains why every formula has by definition a linear negation which expresses a duality or change of
standpoint (Girard 1989):
\bq
``action of type $A$ = reaction of type $A\lood$." (p.77) 
\eq
or in Girard (2000):
\bq
``Concr\`etement la n\'egation correspond \`a la dualit\'e ``action/r\'eaction'' et pas du tout \`a
l'id\'ee de ne pas effectuer une action: typiquement lire/\'ecrire, envoyer/recevoir, sont justicibles de
la n\'egation lin\'eaire.'' (p.532)
\eq
In terms of functional programming or categorical semantics, the negation represents an input-output
duality.  In game terms, it is an opponent-proponent duality.
In Girard (1989) the notion of a reaction of type $A\lood$, as dual to an action of type $A$, is quite
mysterious and not further elaborated.  The only thing Girard mentions is that it
should come down to an ``inversion of causality, i.e. of the sense of time" (Girard 1989).  It is also not
clear to us what this duality would mean in a commutative linear logic context when we consider the case of an action of
type
$\mathfrak{a}
\multimap
\mathfrak{b}$ (or $\mathfrak{b}\lood \multimap \mathfrak{a}\lood$) and a reaction of type $\mathfrak{a} \otimes
\mathfrak{b}\lood$ (or $\mathfrak{b}\lood \otimes \mathfrak{a}$) --- where 
$\mathfrak{a} \multimap \mathfrak{b} = \mathfrak{a}\lood \acht \mathfrak{b}$ and $(\mathfrak{a}\lood \acht
\mathfrak{b})\lood =
\mathfrak{a} \otimes
\mathfrak{b}\lood$.  Thus we tend to agree with 
Girard (1989) where he says that this discussion involves not standard but non-commutative linear logic. 
Indeed, against the background of the causation-propagation duality elaborated upon above, if Girard has
something similar to causation in mind, we know that what is necessary is an implication and
``non-commutative" conjunction which allow us to express causation.  As we will show further on, switching
from standard to non-commutative linear logic effects the meanings of the linear implication and linear
negation.  

\smallskip
Leaving this action-reaction debate aside, we can now
explain why our given interpretation of non-iterability is quite subtle.  First note that in Girard's standard linear
logic, $\vdash \mathfrak{a}$ $\multimap$ $\mathfrak{a}$ is provable from an empty set of premises. 
As such
$\multimap$ represents in $\mathfrak{a} \multimap \mathfrak{a}$ the identity-action which does not really change
resources, but only translates initial ones into final ones.  And although we could perform the action twice in the
following sense:
$\mathfrak{a}_1$ $\multimap$
$\mathfrak{a}_2$ $\multimap$ $\mathfrak{a}_3$ --- for convenience we labeled the resources --- this still
does not count as an iteratable action since $\mathfrak{a}_1$ is to be conceived  as an initial resource,
different from
$\mathfrak{a}_2$, the final resource of the first action.
But we can go further in this discussion and follow Girard in stipulating the fact that we may still encounter situations in
which the picture of ``consuming all initial resources" does not hold. Linear logic henceforth allows also the
expression of those actions which deal with stable situations and which are iteratable. In the latter case the use of
exponentials is necessary where for instance the exponential
$!$ gives $! \mathfrak{a}$ the meaning that
$\mathfrak{a}$'s use as a resource is unlimited.  
\bigskip\par\noindent
In Girard (1989), Girard discusses the
link between states, transitions and the linear implication.  In particular, for us it is the following statement
which places the linear implication in an interesting context, thinking of course about the above discussed
DoQL, (Girard 1989):
\bq
 ``In fact, we would like to represent states by formulas,
and transitions by means of implications of states, in such a way that the state $S'$ is accessible from
$S$ exactly when $S$ $\multimap$ $S'$ is provable from the transitions, taken as
axioms." (p.74)
\eq
In Girard (1989) this statement applies to for instance systems such as Petri nets, Turing
machines, chessboard games, etc.  Focusing ``in this sense" on physical systems,
and using
$\multimap$ for transitions of states, it becomes interesting to investigate how $\multimap$ can be conceived in the context
of our logic of actuality sets.  While a formal comparison on the semantical level will be given in the next
paragraph we have to stress here that there is on the methodological level the following point of difference between
DoQL and linear logic: contrary to DoQL, it should be well understood that linear logic is not
a temporal logic, no preconceptions of time or processes is built into it.  More explicitly
(Girard 1989):
\bq
``Linear logic is eventually about time, space and communication, but is not a temporal logic, or a kind of parallel
language: such approaches try to develop preexisting conceptions about time, processes, etc.  In those
matters, the general understanding is so low that one has good chances to produce systems whose aim is to
avoid the study of their objects [...] The main methodological commitment is to refuse any {\it a priori}
intuition about these objects of study, and to assume that (at least part of) the temporal, the parallel
features of computation are already in Gentzen's approach, but are simply hidden by taxonomy." (p.104)
\eq
As such DoQL started out with a different methodology.  The objects of study are well-known
``scientific objects'' such as physical systems and their properties and the inductions performable on 
physical systems.  In a sense this information has been encrypted in the formulas we used.  All dynamic
propagation- and causation-implications have been labeled by inductions, and this is different from the
linear logical implications which are used to express any  (non-specified) transition. In view of quantum
theory it is indeed the case that one cannot speak about observed quantities without specification of the
particular measurement one performs, and as such, the corresponding induction that encodes von Neumann's
projection postulate, or in more fashionable terms, state-update.  Exactly this could form an argument against
applying the linear logical implications in a context of  physical processes. Thus, we are not tempted to agree with
the proposal in Pratt (1993) of adding linear logical connectives as an extension to quantum logic, but rather focus
on the development of a new logical syntax which will however have some definite similarities with linear
logic, in particular with its quantale semantical fragment.    

\smallskip
In order to get a grasp on the quantale semantics of linear logic we have to say something about its non-commutative
variants. In the literature on non-commutative linear logic, two main directions emerge. In a first direction one introduces
non-commutativity of the multiplicatives by restricting the exchange-rule to circular permutations while in a second
direction one completely drops all structural rules.  Concentrating on the first direction, here linear logic with a cyclic
exchange rule is called {\it cyclic linear logic} and has mainly been developed by D.N Yetter in
Yetter (1990), though we have to note that Girard already makes some remarks on cyclic exchanges in
Girard (1989)\,.  More explicitly we see that the restriction to cyclic permutations means that we consider
the sequents as written on a circle
(Girard 1989)\,.  This then should come down to the fact that $\mathfrak{a}_1 \otimes ... \otimes
\mathfrak{a}_{n-1} \otimes
\mathfrak{a}_n \vdash
\mathfrak{a}_n
\otimes 
\mathfrak{a}_1 \otimes ... \otimes \mathfrak{a}_{n-1}$ is provable in cyclic linear logic.  Of course the meaning of 
$\otimes$ with respect to the standard case changes in the sense that it expresses now ``and then" (Yetter 1990) or when
following Girard (1989) it means that ``in the product $\mathfrak{b} \otimes \mathfrak{a}$, the second
component is done before the first one".   As we will explain in the next paragraph on Girard quantales, it
is exactly the difference between $ -
\otimes
\mathfrak{a}$ and
$\mathfrak{a}
\otimes -$ which leads to the introduction of two different implication-connectives in cyclic linear logic:
$\multimap$ and
$\retro$.  

\smallskip 
Given a unital quantale $({Q}, \bigvee, \otimes)$, with $1$ as the multiplicative neutral 
element with respect to $\otimes$\,, it then follows that 
the endomorphisms $a \otimes -$, $- \otimes a :{Q} \to {Q}$ have right adjoints, $a \multimap - $ and $-
\retro a$ respectively:
$$a \otimes c \leq b\Leftrightarrow c \leq a \multimap b \quad\quad\quad\quad\quad c \otimes a \leq
b\Leftrightarrow c
\leq b
\retro a$$
$$a \multimap b = \bigvee \{ c \in {Q}: a \otimes c \leq b \}\quad\quad b \retro a = \bigvee \{ c \in {Q}: c
\otimes a \leq b \}\,.$$ We know that in the standard linear logic as presented by Girard, the following holds
$\mathfrak{a}
\multimap \mathfrak{b} = 
\mathfrak{a}\lood
\acht
\mathfrak{b} = (\mathfrak{a} \otimes \mathfrak{b}\lood)\lood = (\mathfrak{b}\lood \otimes \mathfrak{a})\lood = \mathfrak{b} \acht \mathfrak{a}\lood =
\mathfrak{b}\lood \multimap \mathfrak{a}\lood$.  In cyclic linear logic where we have now two implications, things change in the
sense that we obtain:
$$
(\mathfrak{a} \otimes \mathfrak{b})\lood = \mathfrak{b}\lood \acht 
\mathfrak{a}\lood\ \ \ \quad\quad  \quad\quad \ \ \ (\mathfrak{a} 
\acht
\mathfrak{b})\lood =
\mathfrak{b}\lood \otimes \mathfrak{a}\lood 
$$
$$\mathfrak{a} \multimap \mathfrak{b} = \mathfrak{a}\lood \acht \mathfrak{b}\ \ \ \quad\quad\quad\quad \quad\quad \ \
\
\mathfrak{b} \retro
\mathfrak{a} = \mathfrak{b} \acht \mathfrak{a}\lood
$$
To be more explicit, the linear negation can be interpreted in the unital quantale $({ Q},\bigvee,\otimes,
1)$ by means of a {\it cyclic dualizing element}, which can be defined as
follows (Rosenthal 1990, Yetter 1990)\,:
\bit
\item
An element $\perp\ \in Q$ is dualizing iff
$ \perp \retro (a \multimap \perp) = a = (\perp \retro a) \multimap \perp$ for all $a\in 
Q$\,. It is cyclic iff
$ a \multimap \perp = \perp \retro a$, for every $a \in {Q}$\,.
\eit
Here the operation $ - \multimap \perp$ or equivalently $ \perp \retro -$ is called the linear negation and can be
written as
$(-)\lood$. Note that a unital quantale with a cyclic dualizing element $\perp$ is called a {\it Girard
quantale}\,, a notion having been introduced in
Yetter (1990)\,.
These Girard quantales can be equipped with modal operators to interpret the linear logical exponentials and form 
as such a
straightforward semantics for the cyclic as well as standard linear logical syntax.  In the latter case $a
\multimap b = b \retro a$.  The disadvantage suffered by cyclic linear logic is that it is still not
``non-commutative enough to properly express time's arrow" (Yetter 1990).  Indeed, if $\multimap$ is to be
conceived as a causal implication then it would have been nice to conceive $\retro$ as expressing past
causality, though this interpretation is too misleading according to (Girard 1989).  In a way we agree with
him since in cyclic linear logic
$\mathfrak{a}
\multimap
\mathfrak{b}$ and $\mathfrak{a}\lood \retro \mathfrak{b}\lood$ are the same --- in the sense that they are both equal
to $\mathfrak{a}\lood \acht \mathfrak{b}$.

\smallskip  
Focusing on the second direction in the non-commutative linear logical literature, we first have to mention
the work of J. Lambek.  Lambek's {\it syntactic calculus} originated in Lambek (1958) and as Girard
admits, is the non-commutative ancestor of linear logic.  However it has to be mentioned that Lambek's syntactic
calculus, as originally developed against a linguistic background, is essentially 
multiplicative and intuitionistic.  Later
on Lambek extended his syntactic calculus with additives and recently renamed his formal calculus {\it bilinear logic}.  In
the same direction we can place the work of V.M. Abrusci who developed a non-commutative version of the 
intuitionistic
linear propositional logic in
Abrusci (1990) and of the classical linear propositional logic in Abrusci (1991).   Specific
to Abrusci's work, however, is the fact that a full removal of the exchange rule requires the introduction of two
different negations and two different implications.  To explain this in detail we switch to the
semantical level of quantales; where we want to note that Abrusci works in Abrusci (1990, 1991) with the more
specific structure of phase spaces which as proved in Rosenthal (1990) are examples of quantales.  This then
leads to the fact that for $a \in { Q}$: $\perp \retro (a \multimap \perp) \not = (\perp \retro a) \multimap a$.  Here
we can follow Abrusci and define $\perp \retro (a \multimap \perp) = {\lood(a\lood)}$ and $ (\perp \retro a) \multimap
\perp = (\lood a)\lood$\,, where on the syntactical level
$\mathfrak{a}\lood$ is called the linear postnegation, $\lood \mathfrak{a}$ the linear retronegation, 
$\mathfrak{a} \multimap
\mathfrak{b}$ the linear postimplication and $\mathfrak{b} \retro \mathfrak{a}$ the linear retroimplication. 
 Further on the
syntactical level we obviously have $\mathfrak{a}\lood \acht \mathfrak{b} = \mathfrak{a} \multimap \mathfrak{b}$ while 
$\mathfrak{b} \acht
\lood\mathfrak{a} = \mathfrak{b} \retro \mathfrak{a}$. 
Let us finish of this paragraph with a note on the fact that there is of course much more to say about (non)-commutative
linear logic, indeed contemporary research is in
full development and heads in the direction of combining cyclic linear logic with
commutative linear logic, we however limit ourselves for the time being to
the overview given.  

\smallskip 
Given the above expositions, it follows that $(DI({\cal L}), \bigvee_{DI({\cal L})}, \otimes_e)$\,, the
quantale fragment emerging for propagation for a specific induction $e$\,, provides an example of the
quantale obtained for commutative linear logic.  Note that the difference of course lies in the fact that in
commutative linear logic no retro-implication different from $\multimap$ is present, not even as an
additional structure.  In this respect, to obtain a retro-implication in linear logic it is necessary to move to
a non-commutative linear logic providing a single quantale in which to interpret $\multimap$ and
$\retro$\,, in sharp contrast to our constructions allowing the interpretation of
$\stackrel{e}{\rightarrow}$ and
$\stackrel{e}{\leftarrow}$ for each specific $e$\,.  
Finally, let us note that while the implication $\multimap$, when focusing on it
as an implication expressing simultaneous consumption of initial resources and production of final resources, is much
stronger than
$\stackrel{e}{\rightarrow}$, it can nevertheless be reconstructed within the framework of DoQL, for which we
refer to Smets (2001).

\bigskip\noindent   
{\bf 5. CONCLUSION}  

\medskip\noindent 
It seems to us that an actual attitude towards the logic of dynamics ought to be pluralistic, as it follows from our
two main paradigmatic examples, dynamic operational quantum logic and Girard's linear logic.  

\smallskip 
The mentioned attempts
that aim to integrate the logic of dynamics as it emerges from for example physical and proof theoretic
considerations fail either on formal grounds or due to conceptual inconsistency. We indeed provided a non-classical
physical counterexample to van Benthem's general dynamic logic, and argued that even for classical systems, from a
physical perspective adjoints rather than relational inverses should be the formal bases for the logic of dynamics.  

\smallskip 
We end by mentioning two promising recent alternative approaches in relating quantum features and linear logic, in
Blute et al (2001) in terms of polycategories and deduction systems
and in Abramsky and C\mbox{oe}cke (2002) in terms of geometry of interaction in categorical format, that is, in terms
of traced monoidal categories.
It is however to soon to obtain conclusions from these lines of thought.

\bigskip\noindent   
{\bf 6. NOTES}

\smallskip\noindent 
{\small
\ben
\itemsep=-4pt  
\itemindent=-0pt            
\labelsep=2pt  
\leftskip=-10pt
\item
From now on OQL will only refer to Geneva
School operational quantum logic.
\item
It is exactly the particular operational foundation of ontological concepts that has caused a
lot of confusion with respect to this approach, including some attacks on it due to
misunderstandings stemming from identification of ``what is", ``what is observed", ``what
will be observed'', ``what would be observed'', ``what could be observed'', etc.  We don't
refer to these papers but cite one that refutes them in a more than convincing way, namely
Foulis and Randall (1984)\,.  We recall here that it were D.J. Foulis and C.H. Randall who
developed an empirical counterpart to C. Piron's ontological approach (Foulis and Randall
1972, Foulis et al 1983, Randall and Foulis 1983). We also quote the review of R. Piziak in
Mathematical Reviews of one of these attacks of Piron's approach exposing their
somewhat doubtful aims (MR86i:81012): 
\begin{quote}``... in fact, they confused the very essence of Piron's system of questions and
propositions, the sharp distinction between properties of a
physical system and operationally testable propositions about the system. [...] In their
reply to Foulis and Randall, HT ignore the list of mathematical errors, confusions and
blunders in their papers except for one (a minor one at that). HT simply reissue their
challenge and dismiss the work of Foulis and Randall as well as Piron as being ``useless from
the physicist's point of view''. However, when one finds an argument in HTM (the main theorem
of their 1981 paper) to the effect that the failure to prove the negation of a theorem
constitutes a proof of the theorem, one might form a different opinion as to whose work is
`useless'. It is right and proper for any scientific work to be scrutinized and criticized
according to its merits. Indeed, this is a main impetus to progress. But if mathematics is
to be used as a tool of criticism, let it be used properly.''\end{quote} See also Smets
(2001) for an overview of most of the criticism and its refutal on sOQL. We refer to Coecke (nd)
for a formulation of sOQL where a conceptually somewhat less rigid, but more general
perspective is proposed, avoiding the notion of test or definite experimental project in the
definition of a physical property as an ontological quality of a system. One of the
motivations for this reformulation is exactly the confusion that the current formulation
seems to cause --- although there is definitely nothing wrong with it as the truly careful
reader knows, on the contrary in fact.
\item
We are not
implying that our scientific theories are to be based on obtained measurement-results.  The
oft-drawn conclusion from this stating that ``ontological existence is independent from any
measurement or observation" also holds in our view.  It is however specific for our
position, which may not be share by every scientific realist, that our {\it knowledge} of
what exists is linked to what we could measure, stated counterfactually. As such we adopt an
``endo perspective'', measurements are not a priori part of our universe of discourse but
incorporated in a conditional way, it is in this sense that e.g. two not simultaneously
observable properties which ontologically exist, can both unproblematically be incorporated
in our description. We refer to Coecke and Smets (2001) for more details on the ``endo versus
exo perspective''.  
A somewhat related view we want to draw to your attention is put forward in Ghins (2000), where he
argues that from the affirmation of some ``existence", under a criterion of existence based on the conditions of
``presence" and ``invariance", certain counterfactuals should reasonably also be affirmed.
\item
We are well aware of the fact that several correspondence theories of truth have been put forward and have also
been criticised.  Adhering to some form of scientific realism does not necessarily imply that one accepts a
correspondence theory, even more it has been suggested that the
debate on the notion of truth can be cut loose from the debate on realism versus anti-realism --- see for instance
Horwich (1997), Tarski (1944).  Still, against the background of sOQL, we are sympathetic towards a contemporary
account of a Tarskian-style semantic correspondence theory --- see for example Niiniluoto (1999, \S3.4).
\item
In
the line of C. Piron (1981) we note that in general, the actual performance of a
definite experimental project can at most serve to prove the falsity of a physicists¹
assumptions, it cannot hand out a prove for them to be true.
\item
Note here the similarity with the generation of a  
quantale structure within the context of process semantics for computational systems sensu Abramsky and
Vickers (1993) and Resende (2000)\,.
\item
Different approaches however do exist for considering potentially destructive
measurements, see for example Faure et al (1995), Amira et al (1998), Coecke and Stubbe (1999), Coecke et
al (2001) and in particular Sourbron (2000)\,.
\item
From this point on we will identify those inductions which have the same action on properties, i.e., we
abstract from the physical procedure to its transitional effect.
\item
A pair of maps $f^*:L\to M$ and $f_*:M\to L$ between posets $L$ and $M$ are {\it Galois adjoint},  
denoted by $f^*\dashv f_*$\,, if
and only if $f^*(a)\leq b\Leftrightarrow a\leq f_*(b)$ if and only if $\forall a\in M:f^*(f_\ast(a))
\leq a$ and $\forall a\in L:a \leq f_{\ast}(f^*(a))$\,.  One could somewhat abusively say that Galois
adjointness generalizes the notion of inverse maps to non-isomorphic objects: in the case that $f^*$ and
$f_*$ are inverse, and thus $L$ and $M$ isomorphic, the above inequalities saturate in equalities.
Whenever
$f^*\dashv f_*$, $f^*$ preserves all existing joins and $f_*$ all existing meets.  This means that for a
Galois adjoint pair between complete lattices, one of these maps preserves all meets and the other
preserves all joins.  
Conversely, for $L$ and
$M$ complete lattices, any meet preserving map $f_*:M\to L$
has a unique join preserving left Galois adjoint
$f^*:a\mapsto\bigwedge\{b\in M|a\leq f_*(b)\}$
and any join preserving map $f^*:L\to M$ a unique meet preserving {right adjoint}
$f_*:b\mapsto\bigvee\{a\in L|f^*(a)\leq b\}$\,.
\item
This dual representation is included in the duality between the
categories \underline{Inf} of complete lattices and meet-preserving maps and \underline{Sup} of complete
lattices and join preserving maps, since this duality is exactly established in terms of Galois adjunction
at the morphism level --- see Coecke et al (2001) for details.
\item 
See also Amira et al (1998) and Coecke and Stubbe (1999) for a
similar development in terms of a so-called operational resolution on the state set.
\item
Let us briefly describe the more rigid argumentation.  Consider the following definitions
for $A\subseteq{\cal L}$\,: i.  
$\bigvee A$ is called disjunctive iff ($\bigvee A$ actual $\Leftrightarrow$ $\exists a\in A: a$ actual)\,;  
ii.  Superposition states for $\bigvee A$ are
 states for which $\bigvee A$ is actual while no $a\in A$ is actual\,;
iii. Superposition properties for $\bigvee A$ are properties $c<\bigvee A$ whose
actuality doesn't imply that at least one $a\in A$ is actual.
We then have that ($\bigvee A$ disjunctive $\Leftrightarrow$ $\bigvee A$
distributive) provided that existence of superposition 
states implies existence of superposition properties (Coecke 2002)\,.
Moreover, any complete lattice $L$ has 
the complete Heyting algebra ${DI}(L)$ of distributive ideals as its distributive hull (Bruns
and Lakser 1970), providing it with a universal property\,. The inclusion preserves all meets and existing
distributive joins. Thus, ${DI}({\cal L})$ encodes all
possible disjunctions of properties, and moreover, it turns out that all ${DI}({\cal L})$-meets are
conjunctive and all
${DI}({\cal L})$-joins are disjunctive --- note that this is definitely not the case in the powerset $P({\cal L})$
of a property lattice, nor in the order ideals $I({\cal L})$ ordered by
inclusion\,.  It follows from this that the object equivalence between: 
i. complete lattices, and, ii. complete Heyting algebras equipped with a
distributive closure (i.e., it preserves distributive sets), encodes an
intuitionistic representation for operational quantum logic --- see Coecke (2002) for details.  
\item
The following map $\hat{e}^*:P(\Sigma)_1\to P(\Sigma)_2$ provides a
counterexample: Let $\Sigma_1=\Sigma_2:=\{p,q,r,s\}$ with as closed subsets
$\bigl\{\emptyset,\{p\},\{q\},\{r\},\{s\},\{q,r,s\},\Sigma\bigr\}\subset P(\Sigma)_1=P(\Sigma)_2$\,, and set
$\hat{e}^*(\{p\})\mapsto\{p,q\}$\,, $\hat{e}^*(\{q\})\mapsto\{q\}$\,, $\hat{e}^*(\{r\})\mapsto\{r\}$\,,
$\hat{e}^*(\{s\})\mapsto\{s\}$\,; one verifies that although $\bigvee\{q,r\}=\bigvee\{r,s\}$ we have
$\bigvee f_{R_e^{-1}}(\{q,r\})=\bigvee\{p,q,r\}\not=\bigvee\{r,s\}=\bigvee f_{R_e^{-1}}(\{r,s\})$ since
$\bigvee\{p,q,r\}$ yields the top element of the property lattice and $\bigvee\{r,s\}$ doesn't\,.
\item
Note here that in
accordance to the literature on linear logic, contra section 2 and section 3 of this paper, we do use the
terms conjunction and disjunction beyond their strict intuitionistic significance.  This however should not
cause any confusion.
\een 
}
\par\vspace{-2.5mm}\noindent  
{\bf 7. ACKNOWLEDGEMENTS} 

\medskip\noindent 
We thank David Foulis, Jim Lambek, Constantin Piron and Isar Stubbe for discussions and comments that have
led to the present content and form of the presentation in this paper.  
Part of the research reflected in this paper was performed by Bob Coecke at
{\sl McGill University, Department of  Mathematics and Statistics, Montreal} and {\sl Imperial College of Science,
Technology \& Medicine, Theoretical Physics Group, London}\,.  Bob Coecke is Postdoctoral Researcher at the 
European TMR Network ``Linear Logic in Computer Science''. Sonja Smets is Postdoctoral Researcher at
Flanders' Fund for Scientific Research. 

\end{document}